\documentclass[12pt]{article}
\input epsf
\usepackage{epsfig}
\usepackage{pict2e}
\usepackage{color}
\usepackage[percent]{overpic}
\usepackage{amsmath}
\usepackage{amsthm}
\usepackage{amssymb}
\usepackage[percent]{overpic}
\def\Ima{\mathrm{Im}\,}
\def\C{\mathbf{C}}
\title{Interactions between Function Theory and Holomorphic Dynamics}
\author{Alexandre Eremenko\thanks{Supported by NSF grant DMS-1665115.}}
\begin{document}
\maketitle

\begin{center}
{\em Dedicated to Walter Bergweiler on the occasion of his 60-th birthday}
\end{center}

\vspace{.1in}
It is not surprising that in the study of dynamics of holomorphic
functions one uses results and tools from the general theory of functions.
Fatou used normal families and distortion theorems.
Douady, Hubbard and Sullivan
pioneered application of quasiconformal mappings. Baker used Ahlfors'
theory to prove that repelling cycles are dense on the Julia set.
Almost all advanced tools of Function Theory (harmonic measure,
hyperbolic metric, extremal length, Nevanlinna theory, Iversen's classification,
refined minimum
modulus estimates for entire functions, subharmonic functions, etc.)
found applications in holomorphic dynamics. 

However there are different, less common kinds of interaction
on which I want to
concentrate in this talk.

1. Development of new tools (and expansion of old ones) in Function Theory
stimulated by applications in Dynamics.

2. Application of dynamical arguments to ``pure'' Function Theory.
I mean here questions of Function Theory which have no dynamic contents
in their formulation. Sometimes it happens that dynamical considerations
help in solving them.

3. Questions of the type mentioned in paragraph 2 may lead to new problems
in dynamics itself.

I will discuss several interesting instances of 1, 2, 3, illustrated by
the work of Walter.
\vspace{.1in}

\noindent
{\bf Generalization of the Wiman--Valiron theory}
\vspace{.1in}

This theory originates from two papers of Wiman of 1914 and 1916,
and during the 20-th century it was one of the main tools of the study of
general entire functions. The original theory, as proposed by Wiman,
and put into more precise form by Valiron,
is stated in terms of an everywhere convergent power series
$$f(z)=\sum_{k=0}^\infty a_nz^n.$$
As the series converges, it has a term of maximal modulus,
this term depends on $|z|$. We define
$$\mu(r,f)=\max\{|a_k|r^k:1\leq k<\infty\}=|a_n|r^n,$$
where $n=n(r,f)$. If there are several $n$ satisfying this equality,
we take the largest of them as $n(r,f)$.

These functions $\mu(r,f)$ and $n(r,f)$ are called the {\em maximal term}
and the {\em central index}, respectively.
The {\em maximum modulus} and a point $z_r$ where it is achieved
are defined by
$$M(r,f)=\max_{|z|\leq r}|f(z)|=|f(z_r)|.$$
It follows from the Maximum Principle that $|z_r|=r$.

The main result of the theory states that for every $\epsilon>0$,
there are arbitrarily large
values of $r$ for which the asymptotic formula   
\begin{equation}\label{main}
f(z)=(1+o(1))\left(\frac{z}{z_r}\right)^{n(r,f)}f(z_r),\quad r\to+\infty.
\end{equation}
holds in the {\em Wiman-Valiron disk}
\begin{equation}\label{disk}
|z-z_r|\leq\frac{r}{(n(r,f))^{1/2+\epsilon}}.
\end{equation}
This formula can be differentiated any number of times:
\begin{equation}\label{diff}
f^{(m)}(z)=(1+o(1))\left(\frac{n(r,f)}{z}\right)^m
\left(\frac{z}{z_r}\right)^{n(r,f)}f(z_r),
\end{equation}
for $z$ in the Wiman-Valiron disk (\ref{disk}).

In fact, these relations hold for most values of $r$: the
{\em exceptional set}
$E\subset(0,+\infty)$ consisting of $r$
for which they do not hold
has finite logarithmic measure, that is
\begin{equation}\label{ex}
\int_E\frac{dr}{r}<\infty.
\end{equation}
These relations show that we have very good control of $f$ in
a disk around the point $z_r$ where the maximum modulus is attained,
and this happens for most $r$. The function $f$ behaves like a 
{\em monomial} in this disk.

One of the main applications of these relations was to the study of
entire solutions of differential equations. When one inserts 
expressions (\ref{diff}) into an algebraic differential equation,
one obtains an asymptotic algebraic equation, from which the growth
rate of $n(r,f)$ can be derived, and thus the growth rate of the
function itself.

Applications of (\ref{main}) to holomorphic dynamics are based on
the observation that the size of the Wiman-Valiron disk
(in the $\log z$-plane)
is $(n(r,f))^{-1/2-\epsilon}$
while the degree of the monomial in the RHS of (\ref{main}) is $n(r,f)$.
This means that the image of the Wiman-Valiron disk
under the monomial (and thus under the function $f$ as well)
completely covers a large ring of the
form $r_1<|z|<r_2$ with $r_1>r$ and $r_2/r_1\to\infty$.
This ring must contain non-exceptional values of $r$ and another
Wiman-Valiron disk. In this way one proves that the {\em escaping set}
$$I(r,f)=\{ z:f^{n}(z)\to\infty\}$$
is non-empty for all entire functions $f$
of degree at least $2$. (This simple argument was 
found independently by Walter and myself,
but I published it earlier.)

The proofs of Wiman and Valiron were substantially simplified my Macintyre
in 1938 who noticed that the function 
\begin{equation}\label{a}
a(r,f):=r\frac{M'(r,f)}{M(r,f)}
\end{equation}
is approximately equal to $$n(r,f)=r\frac{\mu'(r,f)}{\mu(r,f)},$$
and $a(r,f)$ can replace $n(r,f)$ in the main
relation (\ref{main}). This main relation now takes the form
\begin{equation}\label{mc}
f(z)=(1+o(1))\left(\frac{z}{z_r}\right)^{a(r,f)}f(z_r),
\quad\mbox{for}\quad|z-z_r|\leq\frac{r}{(a(r,f))^{1/2+\epsilon}}.
\end{equation}
There is no mentioning of the power series in this statement,
though in
Macintyre's proof it was still used that $f$ is entire.

The ultimate generalization was achieved by Bergweiler, Rippon and Stallard
in 2008. 

\vspace{.1in}
\noindent
{\bf Theorem 1.} {\em Let $D$ be an arbitrary (unbounded) region in the plane,
and $f$ a holomorphic function in $D$ which is unbounded in $D$ but
bounded on $\partial D$. Let
$$M(r,f)=\max\{|f(z)|:z\in D,\;|z|=r\}$$
and let $a(r,f)$ be defined by (\ref{a}).
Then for every $\epsilon>0$ there exists a set $E\subset(0,+\infty)$
satisfying (\ref{ex})
such that for $r\not\in E$ the Wiman--Valiron disk is
contained in $D$
and relation (\ref{mc}) holds.}
\vspace{.1in}

When $D=\C$ we obtain the original statement of Macintyre.
I used this theorem in 1982 in the investigation of differential equations,
but thought erroneously that it was known at that time, and had a wrong proof
of it.

The central part of the argument is the proof of the somewhat surprising
fact that the Wiman-Valiron disk (as in (\ref{mc}))
is always contained in $D$
for non-exceptional values of $r$. Once this is established,
Macintyre's proof applies. This main fact is proved by potential theory,
but we remark that the statement is not true
for arbitrary harmonic
function unbounded in an arbitrary domain and bounded on
the boundary.

The theorem was used to extend some dynamical properties
of entire
functions to meromorphic functions with a direct tract over infinity.
For example: {\em Let $f$ be a meromorphic function with a direct
tract over $\infty$. Then $I(f)$ contains an unbounded component.}
%See the next section for more about direct tracts.

It is remarkable that a question from holomorphic dynamics led to
a dramatic generalization
of one of the main classical tools of the
theory of entire functions.
\vspace{.1in}

\noindent
{\bf Hayman's conjecture on distribution of values of derivatives}
\vspace{.1in}

The conjecture (1967)
was that for every transcendental meromorphic function $f$,
and every $n\geq 1$ the equation
\begin{equation}
\label{haym}
f^n(z)f'(z)=c\quad\mbox{with}\quad c\neq 0
\end{equation}
has infinitely many solutions.
Hayman himself proved this for $n\geq 3$ and Mues for $n=2$. Clunie proved
the conjecture for entire functions. The case
$n=1$ for meromorphic functions remained unsolved for long time.

The paper of Bergweiler and Eremenko (1995) where this was finally proved
(by a unified argument for all $n\geq 1$) uses a variety of diverse tools,
one of them from dynamics. It is interesting that this paper had a lot
of following, generalizing the result in various directions
but no alternative proof was ever found.

The dynamical argument is used in this paper to prove the following statement:
\vspace{.1in}

\noindent
{\bf Theorem 2.}
{\em Suppose that $f$ is a meromorphic function of finite order with infinitely
many multiple zeros.
Then the equation $f'(z)=1$ has infinitely many solutions.}
\vspace{.1in}

{\em Proof.} Consider the meromorphic function $g(z)=z-f(z)$.
It has infinitely many neutral rational fixed points,
namely the multiple zeros
of $f$. To each neutral rational point at least one attracting petal
is attached and the immediate basin of attraction of this petal must contain
either a critical value or an asymptotic value.

If the number of critical points of $g$ is finite,
then the number of critical
values is also finite, and as $g$ has a finite order, the number of asymptotic
values is also finite. In fact the last implication was the main 
result of this paper:
\vspace{.1in}

\noindent
{\bf Theorem 3.} {\em A meromorphic function of finite order with finitely
many critical values has finitely many asymptotic values}.
\vspace{.1in}

As the number of petals is infinite we conclude that the number of critical
points must be infinite, but the critical points of $g$ are solutions
of $f'(z)=1$ which completes the proofs of Theorem 2.
\vspace{.1in}

Recently Walter explained me how to avoid using Fatou's theorem and to prove
a weaker statement then Theorem 2
which is still sufficient for establishing Hayman's
conjecture.
\vspace{.1in}

{\em Let $f$ be a meromorphic function of finite order
with infinitely many zeros, 
all of them multiple. Then the equation $f'(z)=1$ has infinitely many solutions.}
\vspace{.1in}

{\em Proof.} Suppose that this equation has finitely many solutions,
then $g(z):=z-f(z)$ has finitely many critical points. So it has finitely
many asymptotic values by Theorem 3. For such functions, with all finite
critical
and asymptotic values in the disk $|w|\leq R$
we have an inequality of Rippon and Stallard
$$|g(z)|>R^2\quad\Longrightarrow\quad \left|z\frac{g'(z)}{g(z)}\right|\geq
\frac{\log |g(z)|}{16\pi}.$$
This inequality can be traced back to Teichm\"uller.
Any zero $z_k$ of $f$ is a {\em neutral} fixed point of $g$,
$g(z_k)=z_k,\; g'(z_k)=1,$ and $z_k\to\infty$. Inserting these
points to the inequality we obtain a contradiction.
\vspace{.1in}

There were several improvements and generalizations of Hayman's conjecture,
the strongest one is due to Jianming Chang: {\em Let $f$ be a transcendental
meromorphic function whose derivative is bounded on the set of zeros of $f$.
Then the equation $f'(z)=c$ has infinitely many solutions for every $c\neq 0$.}
All generalizations use some forms of theorems 2 and 3 above.

This proof of Hayman's conjecture remains somewhat
mysterious for me. It uses at least four very different ideas. First the
conjecture is reduced to its special case for meromorphic functions
of finite order. This is done with the help of X. Pang's rescaling lemma,
which is a development of Zalcman's lemma. The finite order case
is handled by Theorem 2 which uses some dynamical arguments,
and Theorem 3. Notice that
Theorem 2 itself is not true for functions of infinite
order, as shown by a counterexample. The case when $f$ has finitely many
zeros is treated separately using the deep result of Hayman, whose proof is
based on  completely different ideas. Finally the main technical result,
Theorem 3
uses Iversen's classification of singularities and the argument
from the proof of the Denjoy-Carleman-Ahlfors theorem on direct singularities.

The proof of the Teichm\"uller--Rippon--Stallard
inequality, is the so-called Logarithmic change of the variable,
which is widely used in transcendental dynamics nowadays,
but the statement and proof belong to Function Theory, rather than dynamics.

Theorem 3, according to Mathscinet, Zentralblatt and Google scholar
is the most famous contribution of Walter (and of myself) to Function Theory.
It can be considered a development
of the Denjoy--Carleman--Ahlfors Theorem. 
\vspace{.1in}

\noindent
{\bf Metrics of constant positive curvature with conic singularities}
\vspace{.1in}

Consider the non-linear PDE
\begin{equation}\label{pde}
\Delta u+e^u=8\pi\delta,
\end{equation}
\def\Z{\mathbf{Z}}
where $\delta$ is the sum of $\delta$ functions at the points
of a lattice $\{2m\omega_1+2n\omega_2:m,n\in\Z\}$.
One is interested in doubly periodic solutions, with primitive
periods $2\omega_1,2\omega_2$. As usual, we assume
that $\tau:=\omega_2/\omega_1$ has positive imaginary part.
For which $\tau$ such solutions exist, and how many of them?

Equations of this type are relevant for geometry and physics.
In geometry it describes a conformal metric on the torus, having curvature
$1$ at all points except the origin,
and a conic singularity at the origin with conic angle $3\pi$.
The length element of the metric is $\rho(z)|dz|=2^{-1/2}e^{u/2}|dz|$.

It turns out that a solution exists if and only if $\tau$ satisfies
certain inequalities, and when exists, it comes with a one parametric
family of solutions. This result is originally due to C-S. Lin and Wang,
with a very complicated proof using advanced non-linear PDE theory
and modular forms. I am going to explain a simple proof using
(anti)-meromorphic dynamics.

Locally, near a non-singular point, the general solution is given by
the Liouville formula
\begin{equation}\label{sol}
u(z)=\log\frac{8|f'|^2}{(1+|f|^2)^2},
\end{equation}
where $f$ is any locally univalent meromorphic function. This formula expresses
the fact that a surface with a metric of curvature $1$ is locally
isometric to regions on the standard unit sphere whose metric
has length element $2|dz|/(1+|z|^2).$

This function $f$ has critical points at the lattice points
these critical points are double, and there are no other
critical points. 
The condition that $u$ is periodic translates
in terms of $f$ to the condition
\begin{equation}\label{2nd}
f(z+2\omega_j)=\lambda_jf(z),\quad j=1,2,\quad |\lambda_j|=1.
\end{equation}
So $f$ must be an ``elliptic function of the second kind''
whose both multipliers have absolute value $1$, and which has a single double
critical point at the origin on the torus.
Once we have such a function, formula (\ref{sol}) defines a doubly
periodic solution of (\ref{pde}). It is clear that $f$ can be
multiplied on
any non-zero constant and all relevant properties will
be preserved. So every solution of (\ref{pde}) comes with a one-parametric
family.

The general form of an elliptic function of the second kind with a single
critical point of multiplicity $2$ at the origin is
\begin{equation}\label{fo}
f(z)=ce^{2z\zeta(a)}\frac{\sigma(z-a)}{\sigma(z+a)},
\end{equation}
where $a$ and $c$ are parameters, $c\neq 0$. We use the standard
notation of Weierstrass for $\sigma$, $\zeta$ and $\wp$ functions.

To derive the condition that the multipliers $\lambda$ have
absolute value $1$, we use the formula
$$\sigma(z+2\omega)=-e^{2\eta(z+\omega)}\sigma(z),\quad
\eta=\zeta(\omega),$$
which gives
$$f(z+2\omega)=e^{4\omega\zeta(a)-4\eta a}f(z).$$
So $|\lambda_j|=1,\; j=1,2$ if and only if
both expressions
$$\omega_j\zeta(a)-\eta_ja\quad\mbox{are pure imaginary for}\quad j=1,2.$$
This means that two linear equations with respect to
$a$ and $\zeta=\zeta(a)$ hold:
$$\omega_j\zeta+\overline{\omega_j\zeta}-a\eta_j-\overline{a\eta_j}=0,\quad j=1,2.$$
By eliminating $\overline{\zeta}$ from these two equations we obtain
a single equation of the form
\begin{equation}\label{anti}
Aa+B\overline{a}+\zeta(a)=0
\end{equation}
where
$$A=\frac{\pi}{4\omega_1^2\Ima\tau}-\frac{\eta_1}{\omega_1},\quad
B=-\frac{\pi}{4|\omega_1|^2\Ima\tau}$$
are the unique constants which make the LHS of
(\ref{anti}) invariant under the substitutions
$$(a,\zeta)\mapsto(a+2\omega_j,\zeta+2\eta_j),\quad j=1,2.$$
Equation (\ref{anti}) has to be solved with respect to $a$.

We may notice that this equation has the following potential-theoretic
interpretation. Let $G$ be the ``Green function'' of the torus,
that is a doubly periodic solution of the linear PDE
$$\Delta G=-\delta+c,$$
which exists iff $c=1/(\mbox{area of the torus})$, and $G$ 
is defined up to an additive constant.
The explicit formula 
$$G(z)=-\frac{1}{2\pi}\log|\theta_1(z)|+\frac{(\Ima z)^2}{2\Ima\tau}+C.$$
Solutions $a$ of (\ref{anti}) are exactly the critical points
of $G$. Three of them are
the half-periods $\omega_1,\omega_2,\omega_1+\omega_2$.
which we call the trivial solutions. These trivial solutions,
when inserted into (\ref{fo}) give constant $f$. So we need 
to search for non-trivial solutions.
Here is where dynamics is used.
\vspace{.1in}

\noindent
{\bf Theorem 4.} {\em Equation (\ref{anti}) has  one
pair $(a,-a)$ of non-trivial
solutions if 
\begin{equation}\label{ineq}
\Ima\left(\frac{\pi i}{e_j\omega_1^2+\eta_1\omega_1}-2\tau\right)<0,\quad
j\in\{1,2,3\}.
\end{equation}
and no non-trivial solutions otherwise. 
}
\vspace{.1in}

The region in the $\tau$-halfplane described by (\ref{ineq})
is a curvilinear triangle with zero angles. The image of this region
on the moduli space of tori is bounded by a single analytic curve.
\vspace{.1in}

{\em Sketch of the proof.} We rewrite (\ref{anti}) as a fixed
point condition for an anti-meromorphic function $g$:
$$z=g(z):=-\frac{1}{B}\left(\overline{\zeta(z)}+\overline{Az}\right),$$
satisfying 
$$g(z+2\omega)=g(z)+2\omega.$$
This function has two critical points modulo periods: they are solutions
of $\wp(z)=A$. One can also show that there are no asymptotic values
on the Fatou set. Thus by Fatou's theorem, the number of attracting fixed
points is at most $2$. To estimate the total number of fixed points
we consider the map $\phi(z)=z-g(z)$ as a map of the torus to the Riemann
sphere. The topological degree of this map is $-1$. If $N^+$ and $N^-$
are the numbers of solutions of $\phi(z)=0$ where orientation is
preserved and reversed, then $N^+-N^-=-1$. On the other hand, orientation
preserving solutions are attracting fixed points of $g$, thus $N^+\leq 2$
and $N^-\leq 3$. Thus equation (\ref{anti}) has at most $5$ solutions.
Three of them are trivial, and whenever $a$ is a solution,
$-a$ is also a solution. 

To obtain (\ref{ineq}) we consider two cases. If all three trivial solutions
are repelling fixed points of $g$, then $N^-\geq 3$, which implies that
$N^-=3,\; N^+=2$. If one trivial solution is attracting then it attracts
an orbit of a critical point but then it must attract both critical orbit since the function $g$ is odd. So by Fatou's theorem $N^+=1,\; N^-=2$.

We conclude that (\ref{anti}) has non-trivial solutions if and only
if all three trivial solutions are repelling fixed points.
The set of parameters $\tau$ where this happens is described by 
inequalities (\ref{ineq}).
\vspace{.1in}

\noindent
{\bf Corollary.} {\em Green's function on a torus can have five or
three critical points, depending on whether inequalities (\ref{ineq})
hold or not.}
\vspace{.1in}

It is remarkable that this simply stated fact about an explicit equation
(\ref{anti}) was only discovered in 2010, by Lin and Wang.

The method of the proof, combining Fatou's theorem with topological degree
was introduced by Khavinson and Swiatek who studied the equation
$$z=\overline{f(z)}$$
with a rational function $f$.

Results of this section show that such generalizations of dynamics of
entire functions as dynamics of meromorphic and even anti-meromorphic functions
can be useful in applications.

The family $z\mapsto g(z)$ which is considered here has a remarkable feature
from the point of view of dynamics: the instability locus in the parameter plane
is a simple analytic curve. No fractals, no bifurcations! 
One can make this system holomorphic by replacing $g$ by its second iterate,
but the dependence on parameter $\tau$ will be still only real-analytic.
%Picture
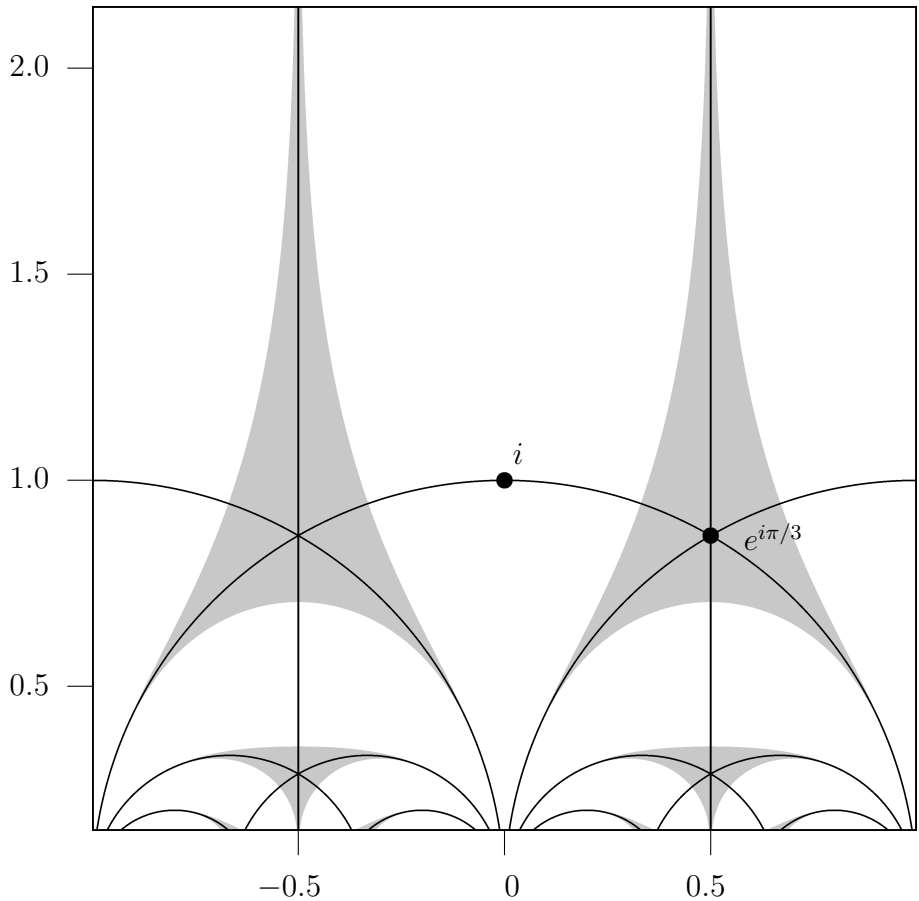
\begin{figure}[htb]
\begin{center}
\begin{overpic}[width=0.80\textwidth]{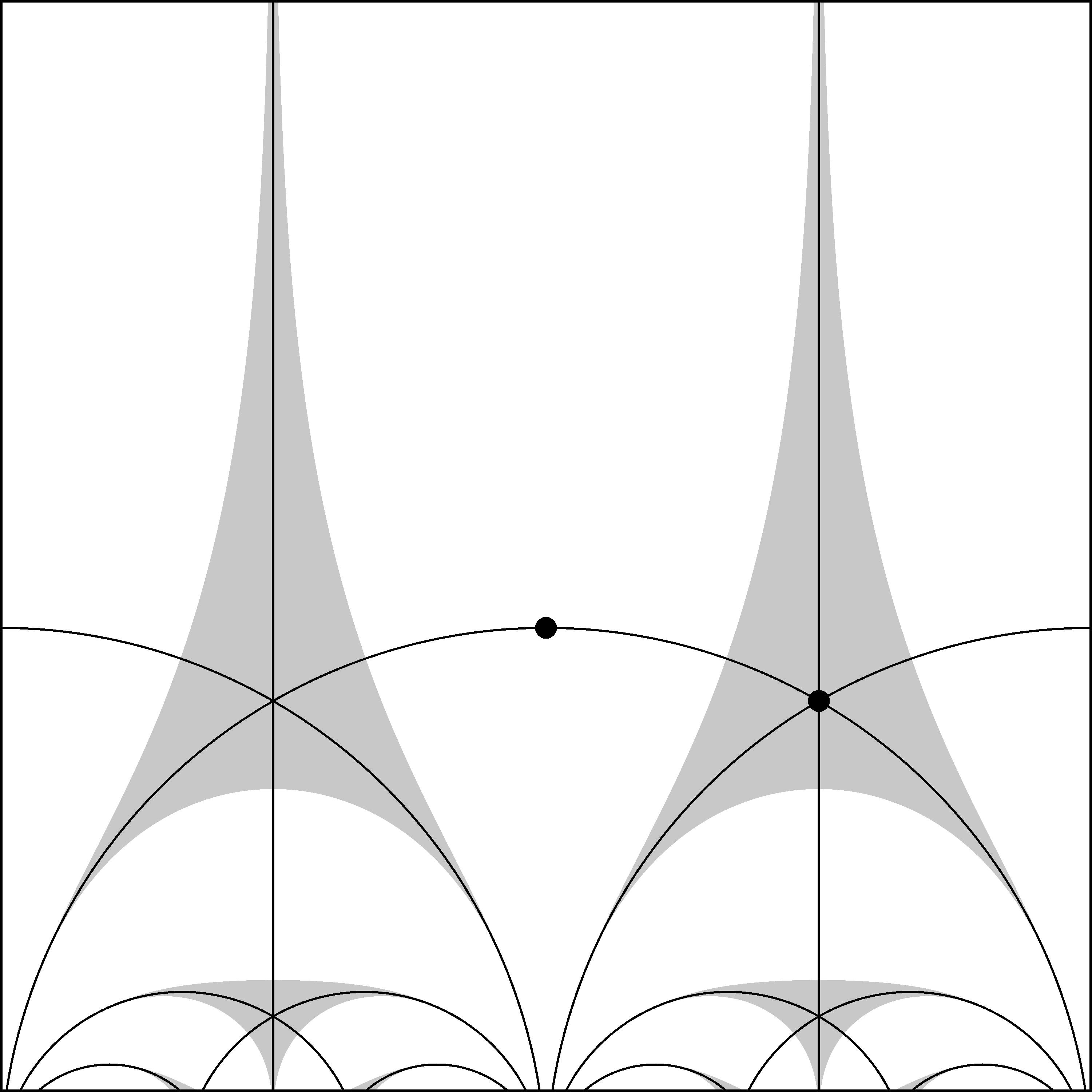}
 \put (25,0) {\line(0,-1){3}}
 \put (50,0) {\line(0,-1){3}}
 \put (75,0) {\line(0,-1){3}}
 \put (20,-8) {$-0.5$}
 \put (50,-8) {$0$}
 \put (72,-8) {$0.5$}
 \put (0,17.5) {\line(-1,0){3}}
 \put (0,42.5) {\line(-1,0){3}}
 \put (0,67.5) {\line(-1,0){3}}
 \put (0,92.5) {\line(-1,0){3}}
 \put (-10,16.5) {$0.5$}
 \put (-10,41.5) {$1.0$}
 \put (-10,66.5) {$1.5$}
 \put (-10,91.5) {$2.0$}
 \put (51,44.5) {$i$}
 \put (79,34.0) {$e^{i\pi/3}$}
\end{overpic}
\vspace*{0.7cm}
\caption{The regions given by (\ref{ineq}) in the $\tau$-plane.}
\label{fig1}
\end{center}
\end{figure}

\end{document}